\documentclass{article}
\usepackage{amssymb}
\usepackage{amsfonts}
\usepackage{amsmath}

\begin{document}

\begin{center}
{\Large Wajsberg algebras of order $n,n\leq 9$}

\begin{equation*}
\end{equation*}

Cristina Flaut, \v{S}\'{a}rka Ho\v{s}kov\'{a}-Mayerov\'{a}, Arsham Borumand
Saeid and

Radu Vasile

\begin{equation*}
\end{equation*}
\end{center}

\textbf{Abstract.}{\small \ In this paper, we describe all finite Wajsberg
algebras of order }$n\leq 9${\small .}

\begin{equation*}
\end{equation*}%

\begin{center}
\bigskip
\end{center}

\textbf{1}. \textbf{Introduction}%
\begin{equation*}
\end{equation*}

Residuated lattices were introduced by Dilworth and Ward, through the papers
[Di; 38], [WD; 39]. A residuated lattic is an algebra $\left( X,\vee ,\wedge
,\odot ,\rightarrow ,0,1\right) $ of type $\left( 2,2,2,2,0,0\right) $ with
an order $\leq $ such that:

- $(X,\vee ,\wedge ,0,1)$ be a bounded lattice;

- $\left( X,\odot ,1\right) $ be a commutative ordered monoid;

- $x\leq y\rightarrow z$ if and only if $y\odot x\leq z$, for all $x,y,z\in
X $, that means $\odot $ and $\rightarrow $ form an orderd pairs ([Pi; 07],
Definition 1.1).

A residuated lattice $X$ is called an MV-algebra if and only if the
following supplimentary condition are satisfied:

\begin{equation*}
\left( x\rightarrow y\right) \rightarrow y=\left( y\rightarrow x\right)
\rightarrow x,\text{ for all }x,y\in X\text{,}
\end{equation*}%
see [Tu; 99], Theorem 2.70.

MV-algebras are introduced by C. C. Chang in [CHA; 58] to give\ new proof
for the completeness of the \L ukasiewicz axioms for infinite valued
propositional logic. These algebras appeared in the specialty literature
under some equivalent names: bounded commutative BCK-algebras or Wajsberg
algebras, ([CT; 96]).

Wajsberg algebras were introduced in 1984, by Font, Rodriguez and Torrens,
through the paper \ [FRT; 84] as an alternative model for the infinite
valued \L ukasiewicz propositional logic.

In [BV; 10] the authors gave an algorithm to find the number of all
non-isomorphic residuated lattices of order $n$, with examples for $n\leq 12$%
.

Knowing connections between residuated lattices, MV-algebras and Wajsbeg
algebras, starting from results obtained in [FV; 19], in this paper we give
Representation Theorem for finite Wajsberg algebras, we give a formula for
the number of nonisomorphic Wajsberg algebras types of order $n$, the total
number of Wajsberg algebras of order $n$ and we describe all finite Wajsberg
algebras of order $n\leq 9$. For $n\in \{1,2,...,9\}$, we obtain the same
numbers of nonisomorphic MV-algebras as in [BV; 10], Table 8, by using other
methods, totally different from their method.

\begin{equation*}
\end{equation*}%
\textbf{\ }

\textbf{2. Preliminaries}%
\begin{equation*}
\end{equation*}

\textbf{Definition 2.1.} ([CHA; 58]) An abelian monoid $\left( X,\theta
,\oplus \right) $ is called \textit{MV-algebra} if and only if we have an
operation $"^{\prime }"$ such that:

i) $(x^{\prime })^{\prime }=x;$

ii) $x\oplus \theta ^{\prime }=\theta ^{\prime };$

iii) $\left( x^{\prime }\oplus y\right) ^{\prime }\oplus y=$ $\left(
y^{\prime }\oplus x\right) ^{\prime }\oplus x,$ for all $x,y\in X.$([Mu;
07]). We denote it by $\left( X,\oplus ,^{\prime },\theta \right) $.\medskip

In an MV-algebra, the following multiplications are also defined:%
\begin{equation*}
x\odot y=\left( x^{\prime }\oplus y^{\prime }\right) ^{\prime },
\end{equation*}%
\begin{equation*}
x\ominus y=x\odot y^{\prime }=\left( x^{\prime }\oplus y\right) ^{\prime }.
\end{equation*}

\textbf{Definition 2.2.}([COM; 00], Definition 4.2.1) An algebra $\left(
W,\circ ,\overline{},1\right) $ of type $\left( 2,1,0\right) ~$is called a 
\textit{Wajsberg algebra (}or\textit{\ W-algebra)} if and only if for every $%
x,y,z\in W,$ we have:

i) $1\circ x=x;$

ii) $\left( x\circ y\right) \circ \left[ \left( y\circ z\right) \circ \left(
x\circ z\right) \right] =1;$

iii) $\left( x\circ y\right) \circ y=\left( y\circ x\right) \circ x;$

iv) $\left( \overline{x}\circ \overline{y}\right) \circ \left( y\circ
x\right) =1.\medskip $

\textbf{Remark 2.3. }([COM; 00], Lemma 4.2.2 and Theorem 4.2.5)

i) If $\left( W,\circ ,\overline{},1\right) $ is a Wajsberg algebra,
defining the following multiplications 
\begin{equation*}
x\odot y=\overline{\left( x\circ \overline{y}\right) }
\end{equation*}%
and 
\begin{equation*}
x\oplus y=\overline{x}\circ y,
\end{equation*}%
for all $x,y\in W$, we obtain that $\left( W,\oplus ,\odot ,\overline{}%
,0,1\right) $ is an MV-algebra.

ii) If $\left( X,\oplus ,\odot ,^{\prime },\theta ,1\right) $ is an
MV-algebra, defining on $X$ the operation%
\begin{equation*}
x\circ y=x^{\prime }\oplus y,
\end{equation*}%
it results that $\left( X,\circ ,^{\prime },1\right) $ is a Wajsberg
algebra.\medskip

\textbf{Definition 2.4.} ([CT; 96]) Let $\left( X,\oplus ,^{\prime },\theta
\right) $ be an MV-algebra. The nonempty subset $I\subseteq X$ is called an 
\textit{ideal} in $X$ if and only if the following conditions are satisfied:

i) $\theta \in I$, where $\theta =\overline{1};$

ii) $x\in I$ and $y\leq x$ implies $y\in I;$

iii) If $x,y\in I$, then $x\oplus y\in I$.\medskip

\textbf{Definition 2.5. }([COM; 00], p. 13) An ideal $P$ of the MV-algebra $%
\left( X,\oplus ,^{\prime },\theta \right) $ is a \textit{prime} ideal in $X$
if and only if for all $x,y\in P$ we have $(x^{\prime }\oplus y)^{\prime
}\in P$ or $(y^{\prime }\oplus x)^{\prime }\in P$.\medskip

\textbf{Definition 2.6. (}[GA; 90], p. 56) Let $\left( W,\circ ,\overline{}%
,1\right) $ be a Wajsberg algebra and $I\subseteq W$ be a nonempty subset. $%
I $ is called an \textit{ideal} in $W$ if and only if the following
conditions are fulfilled:

i) $\theta \in I$, where $\theta =\overline{1};$

ii) $x\in I$ and $y\leq x$ implies $y\in I;$

iii) If $x,y\in I$, then $\overline{x}\circ y\in I$.\medskip

Using connections between MV-algebras and Wajsberg algebras, we give below
the notion of a prime ideal in a Wajsberg algebra.\medskip

\textbf{Definition 2.7.} Let $\left( W,\circ ,\overline{},1\right) $ be a
Wajsberg algebra and $P\subseteq W$ be a nonempty subset. $P$ is called a 
\textit{prime} \textit{ideal} in $W$ if and only if for all $x,y\in P$ we
have $(x\circ y)^{\prime }\in P$ or $(y\circ x)^{\prime }\in P$.\medskip

\textbf{Definition 2.8.} ([CHA; 58]) Let $\left( X,\oplus ,^{\prime },\theta
\right) $ \ be an MV-algebra. The distance function defined on the algebra $%
A $ is:%
\begin{eqnarray*}
d &:&X\times X\rightarrow X,d\left( x,y\right) =\left( x\ominus y\right)
\oplus \left( y\ominus x\right) = \\
&=&(x\odot y^{\prime })\oplus \left( y\odot x^{\prime }\right) = \\
&=&\left( x^{\prime }\oplus y\right) ^{\prime }\oplus \left( y^{\prime
}\oplus x\right) ^{\prime }.
\end{eqnarray*}

\textbf{Definition 2.9. }1)\textbf{\ (}[COM; 00], Proposition 1.2.6) Let $I$
be an ideal in an MV-algebra $X$. We define the following binary relation on 
$X$: 
\begin{equation*}
x\equiv _{I}y\text{ if and only if }d\left( x,y\right) \in I,x,y\in X\text{.}
\end{equation*}

We remark that $I=\{x\in X~/~x\equiv _{I}\theta \}$ and the quotient set $%
X/I $ become an MV-algebra.

2) Let $\left( W,\circ ,\overline{},1\right) $ be a Wajsberg algebra and $I$
be an ideal in an $W$. Using Definition 2.8 and connections between
MV-algebras and Wajsberg algebras, we define the following binary relation
on $W$:%
\begin{equation*}
x\equiv _{I}y\text{ if and only if }\left( x\circ y\right) \circ \left(
y\circ x\right) ^{\prime }\in I\text{.}
\end{equation*}

We remark that $I=\{x\in W~/~x\equiv _{I}\theta \}$ and the quotient set $%
W/I $ become a Wajsberg algebra.\medskip

\textbf{Definition 2.10.} Let $\left( X_{1},\oplus ,^{\prime },\theta
\right) $ and $\left( X_{2},\otimes ,\overleftrightarrow{},0\right) $ be two
MV-algebras. A map $f:X_{1}\rightarrow X_{2}$ is a morphism of Wajsberg
algebras if and only if:

1) $f\left( \theta \right) =0;$

2) $f\left( x\oplus y\right) =f\left( x\right) \otimes f\left( y\right) ;$

3) $f\left( x^{\prime }\right) =\overleftrightarrow{f\left( x\right) }.$

If $f$ is a bijection, therefore the algebras $X\,_{1}$ and $X_{2}$ are
isomorphic. We write this $X\,_{1}\simeq X_{2}.\medskip $

\ \textbf{Definition 2.11.} Let $\left( W_{1},\circ ,\overline{},1\right) $
and $\left( W_{2},\cdot ,^{\prime },1\right) $ be two Wajsberg algebras. A
map $f:W_{1}\rightarrow W_{2}$ is a morphism of Wajsberg algebras if and
only if:

1) $f\left( 0\right) =0;$

2) $f\left( x\circ y\right) =f\left( x\right) \cdot f\left( y\right) ;$

3) $f\left( \overline{x}\right) =\left( f\left( x\right) \right) ^{\prime }.$

If $f$ is a bijection, therefore the algebras $W\,_{1}$ and $W_{2}$ are
isomorphic. We write this $W\,_{1}\simeq W_{2}.\medskip $

\textbf{Definition 2.12.} [FRT; 84] If $\left( W,\circ ,\overline{},1\right) 
$ is a Wajsberg algebra, on $W$ we define the following binary relation 
\begin{equation}
x\leq y~\text{if~and~only~if~}x\circ y=1.  \tag{2.1.}
\end{equation}%
This relation is an order relation, called \textit{the natural order
relation on }$W$\textit{.}

\begin{equation*}
\end{equation*}

\textbf{3. Representation Theorem for finite Wajsberg algebras} 
\begin{equation*}
\end{equation*}

\textbf{Remark 3.1. }([FRT; 84], Theorem 19)

We consider\textbf{\ }$\left( X,\leq \right) $\textbf{\ }a finite totally
ordered set, $X=\{x_{0},x_{1},...,x_{n}\}$, with $x_{0}$ the first and $%
x_{n} $ the last element. With this order relation, the following
multiplication $"\circ "$ are defined on $X$:%
\begin{equation}
\left\{ 
\begin{array}{c}
x_{i}\circ x_{j}=1\text{, if }x_{i}\leq x_{j}; \\ 
x_{i}\circ x_{j}=x_{n-i+j}\text{, otherwise;} \\ 
x_{0}=\theta ,x_{n}=1,x\circ \theta =\overline{x}.%
\end{array}%
\right.  \tag{3.1.}
\end{equation}%
The obtained algebra $\left( X,\circ ,\overline{},1\right) $ is a Wajsberg
algebra and this is the only way to define a Wajsberg algebra structure on a
finite totally ordered set such that the induced order relation on this
algebra is given by relation $(2.1)$, with $\overline{x}_{i}=x_{n-1}$%
.\medskip

\textbf{Definition 3.2. ([}FV; 19\textbf{])} Let $\left( W_{1},\circ ,%
\overline{},\theta \right) $ and $\left( W_{2},\cdot ,^{\prime },1\right) $
be two finite Wajsberg algebras. On the Cartesian product of these algebras, 
$W=W_{1}\times W_{2}$, we define the following multiplication $"\nabla "$,%
\begin{equation}
\left( x_{1},x_{2}\right) \nabla \left( y_{1},y_{2}\right) =\left(
x_{1}\circ y_{1},x_{2}\cdot y_{2}\right) ,  \tag{3.2.}
\end{equation}%
The algebra $\left( W,\nabla ,\rceil ,\mathbf{1}\right) $ is also a Wajsberg
algebra, the complement of the element $\left( x_{1},x_{2}\right) $ is $%
\rceil \left( x_{1},x_{2}\right) =$ $\left( \overline{x}_{1},x_{2}^{\prime
}\right) $ and $\mathbf{1}=\left( \theta ,1\right) $.

If we consider $x=\left( x_{1},x_{2}\right) ,y=\left( y_{1},y_{2}\right) \in
W$, then the order relation on the algebra $\left( W,\nabla ,\rceil ,\mathbf{%
1}\right) $ is defined as follow:%
\begin{equation*}
x\leq _{W}y\text{ if and only if }x_{1}\leq _{W_{1}}y_{1}\text{ and }%
x_{2}\leq _{W_{2}}y_{2}\text{.}
\end{equation*}

\textbf{Proposition 3.3.} ([COM; 00], Theorem 1.3.2) \textit{Let} $X$ 
\textit{be a finite MV-algebra. Therefore,} $X$ \textit{is a direct product
of a family of MV-algebras} $\{X_{i}\}_{i\in \{1,2,...,m\}}$ \textit{if and
only if, there is a family} $\{I_{j}\}_{j\in \{1,2,...,m\}}$ \textit{of
ideals in }$X$\textit{\ such that}

1) $X_{j}\simeq X/I_{j}$, \textit{for all} $j\in \{1,2,...,m\}$;

2) $\underset{j\in \{1,2,...,m\}}{\cap }I_{j}=\{0\}$. $\Box \medskip $

\textbf{Proposition 3.4.} ([CHA; 59], Lemma 1) \textit{Let} $X$ \textit{be
an MV algebra and} $P$ \textit{a prime ideal in} $X$. \textit{Therefore, the
quotient algebra} $X/P$ \textit{is a totally ordered MV-algebra.} $\Box
\medskip $

\textbf{Proposition 3.5.} ([CHA; 59], Lemma 3) \textit{Every MV-algebra is a
direct product of totally ordered MV-algebras.}$\Box \medskip $

Using connections between MV algebras and Wajsberg algebras, we obtain the
following Representation Theorem.\medskip

\textbf{Theorem 3.6}.\ 1) \textit{Let} $W$ \textit{be a finite Wajsberg
algebra. Therefore,} $W$ \textit{is a direct product of a family of Wajsberg
algebras} $\{W_{i}\}_{i\in \{1,2,...,m\}}$ \textit{if and only if, there is
a family} $\{I_{j}\}_{j\in \{1,2,...,m\}}$ \textit{of ideals in} $W$ \textit{%
such that}

i) $W_{j}\simeq W/I_{j}$\textit{, for all} $j\in \{1,2,...,m\}$;

ii) $\underset{j\in \{1,2,...,m\}}{\cap }I_{j}=\{0\}$.

2) \textbf{([}FV; 19\textbf{], Theorem 4.8)} \textit{Each finite Wajsberg
algebra is a direct product of totally ordered Wajsberg algebras.\medskip }

\textbf{Proof.}

1) It results from Proposition 3.1 and Remark 2.3.

2) Indeed, if $\{I_{j}\}_{j\in \{1,2,...,m\}}$ are prime ideals in $W$, from
Proposition 3.2 and Remark 2.3, it results that $W$ is a direct product of
totally ordered Wajsberg algebras. $\Box \medskip $

\textbf{Remark 3.7.}

In \textbf{[}FV; 19\textbf{]} was developed an algorithm which generate all
finite Wajsberg algebras. These algebras are direct product of \ totally
ordered algebras. With Theorem 3.4, we can complete this algorithm:

i) Let $n$ the order of the Wajsberg algebra $\left( W,\circ ,\overline{}%
,1\right) $ and 
\begin{equation*}
n=r_{1}r_{2}...r_{t},r_{i}\in \mathbb{N},1<r_{i}<n,i\in \{1,2,...,t\},
\end{equation*}%
be the decomposition of the number $n$ in factors. Since this decomposition
is not unique, the decompositions with the same terms, but with other order
of them in the product, will we counted one time. The number of all such
decompositions will be denoted with $\pi _{n}$.

ii) Since an MV-algebra is finite if and only if it is isomorphic to a
finite product of totally ordered MV algebras, using connections between
MV-algebras and Wajsberg algebras, we obtain that a Wajsberg algebra is
finite if and only if it is isomorphic to a finite product of totally
ordered Wajsberg algebras. ([HR; 99], Theorem 5.2, p. 43).

ii)(\textbf{[}FV; 19\textbf{], }Theorem 4.8\textbf{)} There are only $\pi
_{n}$ nonismorphic, as ordered sets, Wajsberg algebras with $n$ elements. We
obtain these algebras as a finite product of totally ordered Wajsberg
algebras. We denote them with\textit{\ }$\left( \mathcal{W}_{i}^{n},\nabla
_{i},\rceil _{i},\mathbf{1}_{i},\leq _{i}^{n}\right) $, where $\leq _{i}^{n}$%
is the corresponding order relation on $\mathcal{W}_{i}^{n}$, $i\in
\{1,2,...,\pi _{n}\}$. It is clear that if $n$ is prime, therefore we obtain
only totally ordered Wajsberg algebra.

iii) (\textbf{[}FV; 19\textbf{], }Remark 4.9\textbf{)} We denote with $%
\left( \mathcal{W}_{ij}^{n},\nabla _{ij},\rceil _{ij},\mathbf{1}_{i},\leq
_{ij}\right) $ the Wajsberg algebras isomorphic to $\mathcal{W}_{i}^{n}$,
considered as ordered sets, with $\leq _{ij}^{n}$the corresponding order
relation on $\mathcal{W}_{ij}^{n}$. Let $f_{ij}^{n}:\mathcal{W}%
_{i}^{n}\rightarrow \mathcal{W}_{ij}^{n}$ be such an isomorphism of ordered
sets. The Wajsberg structure on the algebra $\mathcal{W}_{ij}^{n}$ is
defined as follows. For $x,y\in $ $\mathcal{W}_{ij}^{n}$ and $a,b\in 
\mathcal{W}_{i}^{n}$ with $f_{ij}^{n}\left( a\right) =x$ and $%
f_{ij}^{n}\left( b\right) =y$, we define%
\begin{equation*}
x\nabla _{ij}y=f_{ij}^{n}\left( a\right) \nabla _{ij}f_{ij}^{n}\left(
b\right) \overset{def}{=}f_{ij}^{n}\left( a\nabla _{i}b\right) \text{.}
\end{equation*}%
In this way we define a Wajsberg algebra structure on $\mathcal{W}_{ij}^{n}$
such that the induced order relation on this algebra is $\nabla _{ij}$. We
remark that the algebras $\mathcal{W}_{i}^{n}$ and $\mathcal{W}_{ij}^{n}$
are isomorphic as ordered sets and are not always isomorphic as Wajsberg
algebras.

iv) Using Theorem 3.4 and Definition 2.7, a finite Wajsberg algebra is a
direct product of Wajsberg algebras of the form $W/P_{j}$\textit{,} $j\in
\{1,2,...,m\}$ with $\underset{j\in \{1,2,...,m\}}{\cap }P_{j}=\{0\}$ and $%
P\,_{j}$ prime ideals in Wajsberg algebra $W$.\medskip

\textbf{Proposition 3.8.} \textit{Let} $\left( W_{1},\circ ,\overline{}%
,\theta \right) $ \textit{and} $\left( W_{2},\cdot ,^{\prime },1\right) $ 
\textit{be two Wajsberg algebras. If} $f:W_{1}\rightarrow W_{2}$ \textit{is
a bijective map such that} $f\left( \theta \right) =1$ \textit{and} $f\left(
x\circ y\right) =f\left( x\right) \cdot f\left( y\right) $, \textit{therefore%
} $f$ \textit{and} $f^{-1}$ \textit{are morphisms of ordered sets.\medskip }
\ 

\textbf{Proof.} If $x\leq _{W_{1}}y$, therefore $x\circ y=\theta $, we have
that $f\left( x\circ y\right) =f\left( x\right) \cdot f\left( y\right)
=f\left( \theta \right) =1$, therefore $f(x)\leq _{W_{2}}f(y)$. Conversely
is also true.$\Box $ 
\begin{equation*}
\end{equation*}

4. \textbf{Examples of Wajsberg algebras of order }$n\leq 9.$%
\begin{equation*}
\end{equation*}

In this section, we will describe all Wajsberg algebras of order $n\leq 9$.
We will give the number of these algebras, a complete description and a way
to obtain them. For all these algebras, we also compute the prime ideals and
we provide the decompositions given in Theorem 3.6.\bigskip

4.1. \textbf{Wajsberg algebras of order }$4.$\bigskip

\textbf{1)} \textbf{Totally ordered case.} Let $W=\{O\leq A\leq B\leq E\}$
be a totally ordered set. On $W$ we define a multiplication as in relation $%
\left( 3.1\right) $. We have $\overline{A}=B$ and $\overline{B}=A$.
Therefore the algebra $W$ has the following multiplication table:

\begin{equation*}
\begin{tabular}{l|llll}
$\nabla _{0}^{4}$ & $O$ & $A$ & $B$ & $E$ \\ \hline
$O$ & $E$ & $E$ & $E$ & $E$ \\ 
$A$ & $B$ & $E$ & $E$ & $E$ \\ 
$B$ & $A$ & $B$ & $E$ & $E$ \\ 
$E$ & $O$ & $A$ & $B$ & $E$%
\end{tabular}%
\text{.}
\end{equation*}

This algebra has no proper ideals.\medskip

\textbf{2)} \textbf{Partially ordered case.} There is only one type of
partially ordered Wajsberg algebra with $4$ elements, up to an isomorphism
of ordered sets and Wajsberg algebras. Indeed, let $\left(
W_{1}=\{0,1\},\circ ,\overline{},1\right) $ and $\left( W_{2}=\{0,e\},\cdot
,^{\prime },e\right) $ be two finite totally ordered Wajsberg algebras. We
consider $\mathcal{W}_{11}^{4}=W_{1}\times W_{2}=\{\left( 0,0\right) ,\left(
0,e\right) ,\left( 1,0\right) ,\left( 1,e\right) \}=$\newline
$=\{O,A,B,E\}$. On $W_{1}\times W_{2}$ we obtain a Wajsberg algebra
structure by defining the multiplication as in relation $\left( 3.2\right) $%
, (see [FV; 19]). We give this multiplication in the following table:

\begin{equation*}
\begin{tabular}{l|llll}
$\nabla _{11}^{4}$ & $O$ & $A$ & $B$ & $E$ \\ \hline
$O$ & $E$ & $E$ & $E$ & $E$ \\ 
$A$ & $B$ & $E$ & $B$ & $E$ \\ 
$B$ & $A$ & $A$ & $E$ & $E$ \\ 
$E$ & $O$ & $A$ & $B$ & $E$%
\end{tabular}%
\text{.}
\end{equation*}%
\newline

All proper ideals are $P_{1}=\mathbf{\{}O,A\mathbf{\}}$ and $P\,_{2}=\{O,B%
\mathbf{\}}$. These ideals are also prime ideals. We obtain $\mathcal{W}%
_{11}^{4}/P_{1}=\{\overline{O},\overline{E}\}$ and $\mathcal{W}%
_{11}^{4}/P_{2}=\{\overleftrightarrow{O},\overleftrightarrow{E}\}$, where $%
\overline{O}=\{O,A\},\overline{E}=\{B,E\},\overleftrightarrow{O}=\{O,B\},%
\overleftrightarrow{E}=\{A,E\}$. Indeed, since $(O\circ A)\circ (A\circ
O)\prime =E\circ A\prime =E\circ A=A\in P_{1}$, it results $O\equiv
_{P_{1}}A,$ therefore $\overline{O}=\{O,A\}$. The same computations give us
that $(E\circ B)\circ (B\circ E)\prime =B\circ E\prime =B\circ O=A\in P_{1}$
and $E\equiv _{P_{1}}B$. We remark that $(O\circ B)\circ (B\circ O)\prime
=E\circ A\prime =E\circ B=B\notin P_{1}$, etc. From here, we obtain that $%
\mathcal{W}_{11}^{4}=W_{1}\times W_{2}\simeq \mathcal{W}_{11}^{4}/P_{1}%
\times \mathcal{W}_{11}^{4}/P_{2}$, as in Remark 3.7, iv).

If we consider the map 
\begin{eqnarray*}
f_{1} &:&W_{1}\times W_{2}\rightarrow W_{1}\times W_{2} \\
f_{1}\left( A\right) &\text{=}&B,f_{1}\left( B\right) \text{=}A,f_{1}\left(
O\right) \text{=}O,f_{1}\left( E\right) \text{=}E,
\end{eqnarray*}%
we obtain on $W_{1}\times W_{2}$ the same Wajsberg structure, as we can see
in the below table:

\begin{equation*}
\begin{tabular}{l|llll}
$\nabla ^{\prime }$ & $O$ & $A$ & $B$ & $E$ \\ \hline
$O$ & $E$ & $E$ & $E$ & $E$ \\ 
$A$ & $B$ & $E$ & $B$ & $E$ \\ 
$B$ & $A$ & $A$ & $E$ & $E$ \\ 
$E$ & $O$ & $A$ & $B$ & $E$%
\end{tabular}%
\end{equation*}

Therefore, there are only two nonisomorphic Wajsberg algebras of order $4$,
(as W-algebras and as ordered sets). Using connections between MV-algebras
and Wajsberg algebras, we obtain that there are two nonisomorphic
MV-algebras of order $4$. The same number was found in [BV; 10], Table 8, by
using other method.$\bigskip $

4.2. \textbf{Wajsberg algebras of order }$6\bigskip $

1) \textbf{Totally ordered case. }Let $W=\{O\leq A\leq B\leq C\leq D\leq E\}$
be a totally ordered set. On $W$ we define a multiplication as in relation $%
\left( 3\mathbf{.}1\right) $. We have $\overline{A}=D$, $\overline{B}=C$, $%
\overline{C}=B$, $\overline{D}=A$. Therefore the algebra $W$ has the
following multiplication table:%
\begin{equation*}
\begin{tabular}{l|llllll}
$\nabla _{0}^{6}$ & $O$ & $A$ & $B$ & $C$ & $D$ & $E$ \\ \hline
$O$ & $E$ & $E$ & $E$ & $E$ & $E$ & $E$ \\ 
$A$ & $D$ & $E$ & $E$ & $E$ & $E$ & $E$ \\ 
$B$ & $C$ & $D$ & $E$ & $E$ & $E$ & $E$ \\ 
$C$ & $B$ & $C$ & $D$ & $E$ & $E$ & $E$ \\ 
$D$ & $A$ & $B$ & $C$ & $D$ & $E$ & $E$ \\ 
$E$ & $O$ & $A$ & $B$ & $C$ & $D$ & $E$%
\end{tabular}%
.
\end{equation*}

This algebra has no proper ideals.\medskip

\textbf{2)} \textbf{Partially ordered case. }There is only one type of
partially ordered Wajsberg algebra with $6$ elements, up to an isomorphism
of ordered sets. Indeed, $\pi _{6}=1$. Let $\left( W_{1}=\{0,1\},\circ ,%
\overline{},1\right) $ and $\left( W_{2}=\{0,b,e\},\cdot ,^{\prime
},e\right) $ be two finite totally ordered Wajsberg algebras. Using relation 
$\left( 3.1\right) $, on $W_{2}$ we have that $b^{\prime }=b$. We consider $%
W_{1}\times W_{2}=\{\left( 0,0\right) ,\left( 0,b\right) ,\left( 0,e\right)
,\left( 1,0\right) ,\left( 1,b\right) ,\left( 1,e\right) \}=$\newline
$=\{O,A,B,C,D,E\}$. On $W_{1}\times W_{2}$ we obtain a Wajsberg algebra
structure by defining the multiplication as in relation $\left( 3.2\right) $%
. We give this multiplication in the following table:

\begin{equation}
\begin{tabular}{l|llllll}
$\nabla _{11}^{6}$ & $O$ & $A$ & $B$ & $C$ & $D$ & $E$ \\ \hline
$O$ & $E$ & $E$ & $E$ & $E$ & $E$ & $E$ \\ 
$A$ & $D$ & $E$ & $E$ & $D$ & $E$ & $E$ \\ 
$B$ & $C$ & $D$ & $E$ & $C$ & $D$ & $E$ \\ 
$C$ & $B$ & $B$ & $B$ & $E$ & $E$ & $E$ \\ 
$D$ & $A$ & $B$ & $B$ & $D$ & $E$ & $E$ \\ 
$E$ & $O$ & $A$ & $B$ & $C$ & $D$ & $E$%
\end{tabular}%
\ .  \tag{4.1.}
\end{equation}%
\newline
(see relation $\left( 4.4\right) $ from [FV; 19], Example 4.12).

We remark that $A\leq B,A\leq D,C\leq D$ and the other elements can't be
compared in the algebra $\mathcal{W}_{1}^{6}=\mathcal{W}_{11}^{6}=\left(
W_{1}\times W_{2},\nabla _{11}^{6}\right) $. We denote this order relation
with $\leq _{11}^{6}$. All proper ideals are $P_{1}=\{O,A,B\},P_{2}=\{O,C\}$
and are prime ideals. We obtain $\mathcal{W}_{11}^{6}/P_{1}=\{%
\overleftrightarrow{O},\overleftrightarrow{E}\}$ and $\mathcal{W}%
_{11}^{6}/P_{2}=\{\overline{O},\overline{A},\overline{E}\}$, where $%
\overleftrightarrow{O}=\{O,A,B\},\overleftrightarrow{E}=\{C,D,E\},\overline{O%
}=\{O,C\},\overline{A}=\{A,D\},\overline{E}=\{B,E\}$. Indeed, we have $%
\left( E\circ B\right) \circ \left( B\circ E\right) ^{\prime }=B\circ
E^{\prime }=B\circ O=C,$\newline
$\left( A\circ D\right) \circ \left( D\circ A\right) ^{\prime }=E\circ
B^{\prime }=C$ and so on. From here, it results that $\mathcal{W}%
_{11}^{6}\simeq \mathcal{W}_{11}^{6}/P_{1}\times \mathcal{W}_{11}^{6}/P_{2}$%
, as in Remark 3.7, iv).

If we consider the isomorphism 
\begin{eqnarray*}
f_{12}^{6} &:&\left( W_{1}\times W_{2},\nabla _{11}^{6}\right) \rightarrow
\left( W_{1}\times W_{2},\nabla _{12}^{6}\right) , \\
f_{12}^{6}\left( A\right) &\text{=}&A,f_{12}^{6}\left( B\right) \text{=}%
C,f_{12}^{6}\left( C\right) \text{=}B,f_{12}^{6}\left( D\right) \text{=}%
D,f_{12}^{6}\left( O\right) \text{=}O,f_{12}^{6}\left( E\right) \text{=}E,
\end{eqnarray*}%
we obtain on $W_{1}\times W_{2}$ a new Wajsberg algebra structure, with the
multiplication $\nabla _{12}^{6}$ given in relation $\left( 4.5\right) $
from [FV;19], Example 4.12. We denote this algebra with $\mathcal{W}%
_{12}^{6}=\left( W_{1}\times W_{2},\nabla _{12}^{6}\right) $. Using the
above isomorphism, we get that: -algebras $\mathcal{W}_{11}^{6}=\left(
W_{1}\times W_{2},\nabla _{11}^{6}\right) $ and $\mathcal{W}_{12}^{6}=\left(
W_{1}\times W_{2},\nabla _{12}^{6}\right) $ are isomorphic as Wajsberg
algebras;\newline
- $A\leq C,A\leq D,B\leq D$ and the other elements can't be compared in the
algebra $\mathcal{W}_{12}^{6}$. We denote this order relation with $\leq
_{12}^{6}$;\newline
-all proper ideals $P_{1}=\{O,A,C\}$, $P_{2}=\{O,B\}$ are prime ideals;%
\newline
- $\mathcal{W}_{12}^{6}/P_{1}=\{\overleftrightarrow{O},\overleftrightarrow{E}%
\}$ and $\mathcal{W}_{12}^{6}/P_{2}=\{\overline{O},\overline{A},\overline{E}%
\}$, where $\overleftrightarrow{O}=\{O,A,C\},\overleftrightarrow{E}%
=\{B,D,E\},\overline{O}=\{O,B\},\overline{A}=\{A,D\},\overline{E}=\{C,E\}$;%
\newline
-$\mathcal{W}_{12}^{6}\simeq \mathcal{W}_{12}^{6}/P_{1}\times \mathcal{W}%
_{12}^{6}/P_{2}$, as in Remark 3.7, iv).

If we consider the map 
\begin{eqnarray*}
f_{13}^{6} &:&\left( W_{1}\times W_{2},\nabla _{11}^{6}\right) \rightarrow
\left( W_{1}\times W_{2},\nabla _{13}^{6}\right) , \\
f_{13}^{6}\left( A\right) &=&B,f_{13}^{6}\left( B\right) =D,f_{13}^{6}\left(
C\right) =C,f_{13}^{6}\left( D\right) =A,f_{13}^{6}\left( O\right)
=O,f_{13}^{6}\left( E\right) =E,
\end{eqnarray*}%
we obtain on $W_{1}\times W_{2}$ a new Wajsberg algebra structure, with the
multiplication $\nabla _{13}^{6}$ given in relation $\left( 4.6\right) $
from [FV;19], Example 4.12.

We denote this algebra with $\mathcal{W}_{13}^{6}=\left( W_{1}\times
W_{2},\nabla _{13}^{6}\right) $. Using the above map and Proposition 3.8
from above, we get that:\newline
-algebras $\mathcal{W}_{11}^{6}=\left( W_{1}\times W_{2},\nabla
_{11}^{6}\right) $ and $\mathcal{W}_{13}^{6}=\left( W_{1}\times W_{2},\nabla
_{13}^{6}\right) $ are isomorphic only as ordered sets;\newline
-$B\leq A,C\leq A,B\leq D$ and the other elements can't be compared in the
algebra $\mathcal{W}_{13}^{6}$. We denote this order relation with $\leq
_{13}^{6}$.\newline
-all proper ideals $P_{1}=\{O,B,D\}$, $P_{2}=\{O,C\}$ are prime ideals;%
\newline
- $\mathcal{W}_{13}^{6}/P_{1}=\{\overleftrightarrow{O},\overleftrightarrow{E}%
\}$ and $\mathcal{W}_{13}^{6}/P_{2}=\{\overline{O},\overline{B},\overline{E}%
\}$, where $\overleftrightarrow{O}=\{O,B,D\},\overleftrightarrow{E}%
=\{C,A,E\},\overline{O}=\{O,C\},\overline{A}=\{A,D\},\overline{E}=\{C,E\}$;%
\newline
-$\mathcal{W}_{13}^{6}\simeq \mathcal{W}_{13}^{6}/P_{1}\times \mathcal{W}%
_{13}^{6}/P_{2}$, as in Remark 3.7, iv).

If we consider $W_{2}\times W_{1}=\{\left( 0,0\right) ,\left( 0,1\right)
,\left( b,0\right) ,\left( b,1\right) ,\left( e,0\right) ,\left( e,1\right)
\}=$\newline
$=\{O,A,B,C,D,E\}$, on $W_{2}\times W_{1}$ we obtain a Wajsberg algebra
structure by defining the multiplication as in relation $\left( 3.2\right) $%
. We give this multiplication in the following table:

\begin{equation}
\begin{tabular}{l|llllll}
$\nabla _{14}^{6}$ & $O$ & $A$ & $B$ & $C$ & $D$ & $E$ \\ \hline
$O$ & $E$ & $E$ & $E$ & $E$ & $E$ & $E$ \\ 
$A$ & $D$ & $E$ & $D$ & $E$ & $D$ & $E$ \\ 
$B$ & $C$ & $C$ & $E$ & $E$ & $E$ & $E$ \\ 
$C$ & $B$ & $C$ & $D$ & $E$ & $D$ & $E$ \\ 
$D$ & $A$ & $A$ & $C$ & $C$ & $E$ & $E$ \\ 
$E$ & $O$ & $A$ & $B$ & $C$ & $D$ & $E$%
\end{tabular}%
.  \tag{4.2.}
\end{equation}%
\newline
(see relation $\left( 4.7\right) $ from [FV; 19], Example 4.12)

The algebras $\left( W_{1}\times W_{2},\nabla _{11}^{6}\right) $ and $\left(
W_{2}\times W_{1},\nabla _{14}^{6}\right) $ are isomorphic as Wajsberg
algebras, by taking the map 
\begin{eqnarray*}
f_{14}^{6} &:&\left( W_{1}\times W_{2},\nabla _{11}^{6}\right) \rightarrow
\left( W_{2}\times W_{1},\nabla _{14}^{6}\right) , \\
f_{14}^{6}\left( A\right) &=&B,f_{14}^{6}\left( B\right) =D,f_{14}^{6}\left(
C\right) =A,f_{14}^{6}\left( D\right) =C,f_{14}^{6}\left( O\right)
=O,f_{14}^{6}\left( E\right) =E.
\end{eqnarray*}%
In $\mathcal{W}_{14}^{6}=\left( W_{2}\times W_{1},\nabla _{14}^{6}\right) $,
we have $A\leq C,B\leq C,B\leq D~$and the other elements can't be compared.
We denote this order relation with $\leq _{14}^{6}$. All proper ideals $%
P_{1}=\{O,B,D\}$, $P_{2}=\{O,A\}$ are prime ideals.

We also have:\newline
-$\mathcal{W}_{14}^{6}/P_{1}=\{\overleftrightarrow{O},\overleftrightarrow{E}%
\}$ and $\mathcal{W}_{14}^{6}/P_{2}=\{\overline{O},\overline{B},\overline{E}%
\}$, where $\overleftrightarrow{O}=\{O,B,D\},\overleftrightarrow{E}%
=\{C,A,E\},\overline{O}=\{O,A\},\overline{B}=\{B,C\},\overline{E}=\{D,E\}$;%
\newline
-$\mathcal{W}_{13}^{6}\simeq \mathcal{W}_{13}^{6}/P_{1}\times \mathcal{W}%
_{13}^{6}/P_{2}$, as in Remark 3.7, iv).

If we take the map

\begin{eqnarray*}
f_{15}^{6} &:&\left( W_{1}\times W_{2},\nabla _{11}^{6}\right) \rightarrow
\left( W_{2}\times W_{1},\nabla _{15}^{6}\right) , \\
f_{15}^{6}\left( A\right) &=&B,f_{15}^{6}\left( B\right) =A,f_{15}^{6}\left(
C\right) =D,f_{15}^{6}\left( D\right) =C,f_{15}^{6}\left( O\right)
=O,f_{15}^{6}\left( E\right) =E.
\end{eqnarray*}%
\begin{equation*}
\begin{tabular}{l|llllll}
$\nabla _{15}^{6}$ & $O$ & $A$ & $B$ & $C$ & $D$ & $E$ \\ \hline
$O$ & $E$ & $E$ & $E$ & $E$ & $E$ & $E$ \\ 
$A$ & $D$ & $E$ & $C$ & $C$ & $D$ & $E$ \\ 
$B$ & $C$ & $E$ & $E$ & $E$ & $C$ & $E$ \\ 
$C$ & $B$ & $A$ & $A$ & $E$ & $C$ & $E$ \\ 
$D$ & $A$ & $A$ & $A$ & $E$ & $E$ & $E$ \\ 
$E$ & $O$ & $A$ & $B$ & $C$ & $D$ & $E$%
\end{tabular}%
.
\end{equation*}%
we obtain the Wajsberg algebra $\mathcal{W}_{15}^{6}=\left( W_{2}\times
W_{1},\nabla _{15}^{6}\right) $. In this algebra, we have $B\leq C,B\leq
A,D\leq C~$and the other elements can't be compared. We denote this order
relation with $\leq _{15}^{6}$. The algebras $\mathcal{W}_{11}^{6}=\left(
W_{1}\times W_{2},\nabla _{11}^{6}\right) $ and $\mathcal{W}_{15}^{6}=\left(
W_{2}\times W_{1},\nabla _{15}^{6}\right) $ are isomorphic as Wajsberg
algebras. All proper ideals $P_{1}=\{O,A,B\}$, $P_{2}=\{O,D\}$ are prime.

We also have:\newline
-$\mathcal{W}_{15}^{6}/P_{1}=\{\overleftrightarrow{O},\overleftrightarrow{E}%
\}$ and $\mathcal{W}_{15}^{6}/P_{2}=\{\overline{O},\overline{B},\overline{E}%
\}$, where $\overleftrightarrow{O}=\{O,A,B\},\overleftrightarrow{E}%
=\{C,D,E\},\overline{O}=\{O,D\},\overline{B}=\{B,C\},\overline{E}=\{A,E\}$;%
\newline
-$\mathcal{W}_{15}^{6}\simeq \mathcal{W}_{15}^{6}/P_{1}\times \mathcal{W}%
_{15}^{6}/P_{2}$, as in Remark 3.7, iv).

If we consider the map

\begin{eqnarray*}
f_{16}^{6} &:&\left( W_{1}\times W_{2},\nabla _{11}^{6}\right) \rightarrow
\left( W_{2}\times W_{1},\nabla _{16}^{6}\right) , \\
f_{16}^{6}\left( A\right) &=&C,f_{16}^{6}\left( B\right) =D,f_{16}^{6}\left(
C\right) =B,f_{16}^{6}\left( D\right) =A,f_{16}^{6}\left( O\right)
=O,f_{16}^{6}\left( E\right) =E.
\end{eqnarray*}%
\begin{equation*}
\begin{tabular}{l|llllll}
$\nabla _{16}^{6}$ & $O$ & $A$ & $B$ & $C$ & $D$ & $E$ \\ \hline
$O$ & $E$ & $E$ & $E$ & $E$ & $E$ & $E$ \\ 
$A$ & $C$ & $E$ & $A$ & $D$ & $D$ & $E$ \\ 
$B$ & $D$ & $E$ & $E$ & $D$ & $D$ & $E$ \\ 
$C$ & $A$ & $E$ & $A$ & $E$ & $E$ & $E$ \\ 
$D$ & $B$ & $A$ & $B$ & $A$ & $E$ & $E$ \\ 
$E$ & $O$ & $A$ & $B$ & $C$ & $D$ & $E$%
\end{tabular}%
.
\end{equation*}%
we obtain the Wajsberg algebra $\mathcal{W}_{16}^{6}=\left( W_{2}\times
W_{1},\nabla _{16}^{6}\right) $. In this algebra, we have $C\leq D,C\leq
A,B\leq A~$and the other elements can't be compared. We denote this order
relation with $\leq _{16}^{6}$. The algebras $\mathcal{W}_{11}^{6}=\left(
W_{2}\times W_{1},\nabla _{11}^{6}\right) $ and $\mathcal{W}_{16}^{6}=\left(
W_{2}\times W_{1},\nabla _{16}^{6}\right) $ are isomorphic only as ordered
sets. All proper ideals $P_{1}=\{O,C,D\}$, $P_{2}=\{O,B\}$ are prime.

We also have:\newline
-$\mathcal{W}_{16}^{6}/P_{1}=\{\overleftrightarrow{O},\overleftrightarrow{E}%
\}$ and $\mathcal{W}_{16}^{6}/P_{2}=\{\overline{O},\overline{A},\overline{E}%
\}$, where $\overleftrightarrow{O}=\{O,C,D\},\overleftrightarrow{E}%
=\{A,B,E\},\overline{O}=\{O,B\},\overline{A}=\{A,C\},\overline{E}=\{D,E\}$;%
\newline
-$\mathcal{W}_{16}^{6}\simeq \mathcal{W}_{16}^{6}/P_{1}\times \mathcal{W}%
_{16}^{6}/P_{2}$, as in Remark 3.7, iv).

If we take the map

\begin{eqnarray*}
f_{17}^{6} &:&\left( W_{1}\times W_{2},\nabla _{11}^{6}\right) \rightarrow
\left( W_{2}\times W_{1},\nabla _{17}^{6}\right) , \\
f_{17}^{6}\left( A\right) &=&C,f_{17}^{6}\left( B\right) =A,f_{17}^{6}\left(
C\right) =D,f_{17}^{6}\left( D\right) =B,f_{17}^{6}\left( O\right)
=O,f_{17}^{6}\left( E\right) =E.
\end{eqnarray*}%
\begin{equation*}
\begin{tabular}{l|llllll}
$\nabla _{17}^{6}$ & $O$ & $A$ & $B$ & $C$ & $D$ & $E$ \\ \hline
$O$ & $E$ & $E$ & $E$ & $E$ & $E$ & $E$ \\ 
$A$ & $D$ & $E$ & $B$ & $B$ & $D$ & $E$ \\ 
$B$ & $C$ & $A$ & $E$ & $A$ & $B$ & $E$ \\ 
$C$ & $B$ & $E$ & $E$ & $E$ & $B$ & $E$ \\ 
$D$ & $A$ & $A$ & $E$ & $A$ & $E$ & $E$ \\ 
$E$ & $O$ & $A$ & $B$ & $C$ & $D$ & $E$%
\end{tabular}%
.
\end{equation*}%
we obtain the Wajsberg algebra $\mathcal{W}_{17}^{6}=\left( W_{2}\times
W_{1},\nabla _{17}^{6}\right) $. In this algebra, we have $C\leq A,C\leq
B,D\leq B~$and the other elements can't be compared. The algebras $\mathcal{W%
}_{11}^{6}=\left( W_{1}\times W_{2},\nabla _{11}^{6}\right) $ and $\mathcal{W%
}_{17}^{6}=\left( W_{2}\times W_{1},\nabla _{17}^{6}\right) $ are isomorphic
as Wajsberg algebras.We denote this order relation with $\leq _{17}^{6}$.
All proper ideals $P_{1}=\{O,A,C\}$, $P_{2}=\{O,D\}$ are prime ideals.

We also have:\newline
-$\mathcal{W}_{17}^{6}/P_{1}=\{\overleftrightarrow{O},\overleftrightarrow{E}%
\}$ and $\mathcal{W}_{17}^{6}/P_{2}=\{\overline{O},\overline{B},\overline{E}%
\}$, where $\overleftrightarrow{O}=\{O,A,C\},\overleftrightarrow{E}%
=\{B,D,E\},\overline{O}=\{O,D\},\overline{B}=\{B,C\},\overline{E}=\{A,E\}$;%
\newline
-$\mathcal{W}_{17}^{6}\simeq \mathcal{W}_{17}^{6}/P_{1}\times \mathcal{W}%
_{17}^{6}/P_{2}$, as in Remark 3.7, iv).

From the above, we have that there are only two types of nonisomorphic
Wajsberg algebras of order $6$. Using connections between MV-algebras and
Wajsberg algebras, we obtain that there are two nonisomorphic MV-algebras of
order $6$. The same number was found in [BV; 10], Table 8, by using other
method. Since we gave enough examples, now we wonder how we can find all
Wajsberg algebras of order $6$? For this purpose, we must find all
isomorphisms of ordered sets $f:$ $\mathcal{W}_{11}^{6}\rightarrow \mathcal{W%
}_{11}^{6}$, such that $f\left( O\right) =O,f\left( E\right) =E$. Therefore,
there are $4!=24$ isomorphisms and, in turn, $24$ partially ordered Wajsberg
algebras: $8$ are isomorphic with $\mathcal{W}_{11}^{6}$ as Wajsberg
algebras and ordered sets and the next $16$ are isomorphic with $\mathcal{W}%
_{11}^{6}$ only as ordered sets. In total, there are 25 Wajsberg algebras of
order $6$, as we can see in the below table in which are described all
isomorphisms $f_{1j}^{6}:\mathcal{W}_{11}^{6}\rightarrow \mathcal{W}_{1j}^{6}
$, $j\in \{1,...,24\}$, such that $f_{1j}^{6}\left( x\nabla
_{11}^{6}y\right) =f_{1j}^{6}\left( x\right) \nabla
_{1j}^{6}f_{1j}^{6}\left( y\right) $ .

$%
\begin{tabular}[t]{|l|l|l|}
\hline
$\mathcal{W}_{1j}^{6}$ & {\tiny The isomorphism and order relation} & {\tiny %
Isomorpism with }$\mathcal{W}_{11}^{6}$ \\ \hline
$\mathcal{W}_{11}^{6}$ & $f_{11}^{6}\left( A\right) \text{=}%
A,f_{11}^{6}\left( B\right) \text{=}B,f_{11}^{6}\left( C\right) \text{=}%
C,f_{12}^{6}\left( D\right) \text{=}D$ & {\tiny isomorphism of Wajsberg
algebras} \\ \hline
$\mathcal{W}_{12}^{6}$ & $f_{12}^{6}\left( A\right) \text{=}%
A,f_{12}^{6}\left( B\right) \text{=}C,f_{12}^{6}\left( C\right) \text{=}%
B,f_{12}^{6}\left( D\right) \text{=}D$ & {\tiny isomorphism of Wajsberg
algebras} \\ \hline
$\mathcal{W}_{13}^{6}$ & $f_{13}^{6}\left( A\right) $=$B,f_{13}^{6}\left(
B\right) $=$D,f_{13}^{6}\left( C\right) $=$C,f_{13}^{6}\left( D\right) $=$A$
& {\tiny only isomorphism of ordered sets} \\ \hline
$\mathcal{W}_{14}^{6}$ & $f_{14}^{6}\left( A\right) $=$B,f_{14}^{6}\left(
B\right) $=$D,f_{14}^{6}\left( C\right) $=$A,f_{14}^{6}\left( D\right) $=$C$
& {\tiny isomorphism of Wajsberg algebras} \\ \hline
$\mathcal{W}_{15}^{6}$ & $f_{15}^{6}\left( A\right) $=$B,f_{15}^{6}\left(
B\right) $=$A,f_{15}^{6}\left( C\right) $=$D,f_{15}^{6}\left( D\right) $=$C$
& {\tiny isomorphism of Wajsberg algebras} \\ \hline
$\mathcal{W}_{16}^{6}$ & $f_{16}^{6}\left( A\right) $=$C,f_{16}^{6}\left(
B\right) $=$D,f_{16}^{6}\left( C\right) $=$B,f_{16}^{6}\left( D\right) $=$A$
& {\tiny only isomorphism of ordered sets} \\ \hline
$\mathcal{W}_{17}^{6}$ & $f_{17}^{6}\left( A\right) $=$C,f_{17}^{6}\left(
B\right) $=$A,f_{17}^{6}\left( C\right) $=$D,f_{17}^{6}\left( D\right) $=$B$
& {\tiny isomorphism of Wajsberg algebras} \\ \hline
$\mathcal{W}_{18}^{6}$ & $f_{18}^{6}\left( A\right) $=$C,f_{18}^{6}\left(
B\right) $=$B,f_{18}^{6}\left( C\right) $=$D,f_{18}^{6}\left( D\right) $=$A$
& {\tiny only isomorphism of ordered sets} \\ \hline
$\mathcal{W}_{19}^{6}$ & $f_{19}^{6}\left( A\right) $=$C,f_{19}^{6}\left(
B\right) $=$B,f_{19}^{6}\left( C\right) $=$A,f_{19}^{6}\left( D\right) $=$D$
& {\tiny only isomorphism of ordered sets} \\ \hline
$\mathcal{W}_{1,10}^{6}$ & $f_{1\text{,}10}^{6}\left( A\right) $=$C$,$f_{1%
\text{,}10}^{6}\left( B\right) $=$D$,$f_{1\text{,}10}^{6}\left( C\right) $=$%
A $,$f_{1\text{,}10}^{6}\left( D\right) $=$B$ & {\tiny isomorphism of
Wajsberg algebras} \\ \hline
$\mathcal{W}_{1,11}^{6}$ & $f_{1\text{,}11}^{6}\left( A\right) $=$C,f_{1%
\text{,}11}^{6}\left( B\right) $=$A,f_{1\text{,}11}^{6}\left( C\right) $=$%
B,f_{1\text{,}11}^{6}\left( D\right) $=$D$ & {\tiny only isomorphism of
ordered sets} \\ \hline
$\mathcal{W}_{1,12}^{6}$ & $f_{1,12}^{6}\left( A\right) $=$%
B,f_{1,12}^{6}\left( B\right) $=$C,f_{1,12}^{6}\left( C\right) $=$%
D,f_{1,12}^{6}\left( D\right) $=$A$ & {\tiny only isomorphism of ordered sets%
} \\ \hline
$\mathcal{W}_{1,13}^{6}$ & $f_{1,13}^{6}\left( A\right) $=$%
B,f_{1,13}^{6}\left( B\right) $=$C,f_{1,13}^{6}\left( C\right) $=$%
A,f_{1,13}^{6}\left( D\right) $=$D$ & {\tiny only isomorphism of ordered sets%
} \\ \hline
$\mathcal{W}_{1,14}^{6}$ & $f_{1,14}^{6}\left( A\right) $=$%
B,f_{1,14}^{6}\left( B\right) $=$A,f_{1,14}^{6}\left( C\right) $=$%
C,f_{1,14}^{6}\left( D\right) $=$D$ & {\tiny only isomorphism of ordered sets%
} \\ \hline
$\mathcal{W}_{1,15}^{6}$ & $f_{1,15}^{6}\left( A\right) \text{=}%
A,f_{1,15}^{6}\left( B\right) \text{=}C,f_{1,15}^{6}\left( C\right) \text{=}%
D,f_{1,15}^{6}\left( D\right) \text{=}B$ & {\tiny only isomorphism of
ordered sets} \\ \hline
$\mathcal{W}_{1,16}^{6}$ & $f_{1,16}^{6}\left( A\right) \text{=}%
A,f_{1,16}^{6}\left( B\right) \text{=}D,f_{1,16}^{6}\left( C\right) \text{=}%
C,f_{1,16}^{6}\left( D\right) \text{=}B$ & {\tiny only isomorphism of
ordered sets} \\ \hline
$\mathcal{W}_{1,17}^{6}$ & $f_{1,17}^{6}\left( A\right) \text{=}%
A,f_{1,17}^{6}\left( B\right) \text{=}D,f_{1,17}^{6}\left( C\right) \text{=}%
B,f_{1,17}^{6}\left( D\right) \text{=}C$ & {\tiny only isomorphism of
ordered sets} \\ \hline
$\mathcal{W}_{1,18}^{6}$ & $f_{1,18}^{6}\left( A\right) \text{=}%
A,f_{1,18}^{6}\left( B\right) \text{=}B,f_{1,18}^{6}\left( C\right) \text{=}%
D,f_{1,18}^{6}\left( D\right) \text{=}C$ & {\tiny only isomorphism of
ordered sets} \\ \hline
$\mathcal{W}_{1,19}^{6}$ & $f_{1,19}^{6}\left( A\right) $=$%
D,f_{1,19}^{6}\left( B\right) $=$B,f_{1,19}^{6}\left( C\right) $=$%
C,f_{1,19}^{6}\left( D\right) $=$A$ & {\tiny isomorphism of Wajsberg algebras%
} \\ \hline
$\mathcal{W}_{1,20}^{6}$ & $f_{1,20}^{6}\left( A\right) $=$%
D,f_{1,20}^{6}\left( B\right) $=$C,f_{1,20}^{6}\left( C\right) $=$%
B,f_{1,20}^{6}\left( D\right) $=$A$ & {\tiny only isomorphism of ordered sets%
} \\ \hline
$\mathcal{W}_{1,21}^{6}$ & $f_{1,21}^{6}\left( A\right) $=$%
D,f_{1,21}^{6}\left( B\right) $=$A,f_{1,21}^{6}\left( C\right) $=$%
C,f_{1,21}^{6}\left( D\right) $=$B$ & {\tiny only isomorphism of ordered sets%
} \\ \hline
$\mathcal{W}_{1,22}^{6}$ & $f_{1,22}^{6}\left( A\right) $=$%
D,f_{1,22}^{6}\left( B\right) $=$C,f_{1,22}^{6}\left( C\right) $=$%
A,f_{1,22}^{6}\left( D\right) $=$B$ & {\tiny only isomorphism of ordered sets%
} \\ \hline
$\mathcal{W}_{1.23}^{6}$ & $f_{1,23}^{6}\left( A\right) $=$%
D,f_{1,23}^{6}\left( B\right) $=$A,f_{1,23}^{6}\left( C\right) $=$%
B,f_{1,23}^{6}\left( D\right) $=$C$ & {\tiny isomorphism of Wajsberg algebras%
} \\ \hline
$\mathcal{W}_{1,24}^{6}$ & $f_{1,24}^{6}\left( A\right) $=$%
D,f_{1,24}^{6}\left( B\right) $=$B,f_{1,24}^{6}\left( C\right) $=$%
A,f_{1,24}^{6}\left( D\right) $=$C$ & {\tiny only isomorphism of ordered sets%
} \\ \hline
\end{tabular}%
$\newline

\bigskip \textbf{4.3.} \textbf{Wajsberg algebras of order }$8$

1) \textbf{Totally ordered case. }Let $W=\{O\leq X\leq Y\leq Z\leq T\leq
U\leq V\leq E\}$ be a totally ordered set. On $W$ we define a multiplication
as in relation $\left( 3\mathbf{.}1\right) $. We have $\overline{X}=V$, $%
\overline{Y}=U$, $\overline{Z}=T$. Therefore the algebra $W$ has the
following multiplication table:%
\begin{equation*}
\begin{tabular}{l|llllllll}
$\nabla _{0}^{8}$ & $O$ & $X$ & $Y$ & $Z$ & $T$ & $U$ & $V$ & $E$ \\ \hline
$O$ & $E$ & $E$ & $E$ & $E$ & $E$ & $E$ & $E$ & $E$ \\ 
$X$ & $V$ & $E$ & $E$ & $E$ & $E$ & $E$ & $E$ & $E$ \\ 
$Y$ & $U$ & $V$ & $E$ & $E$ & $E$ & $E$ & $E$ & $E$ \\ 
$Z$ & $T$ & $U$ & $V$ & $E$ & $E$ & $E$ & $E$ & $E$ \\ 
$T$ & $Z$ & $T$ & $U$ & $V$ & $E$ & $E$ & $E$ & $E$ \\ 
$U$ & $Y$ & $Z$ & $T$ & $U$ & $V$ & $E$ & $E$ & $E$ \\ 
$V$ & $X$ & $Y$ & $Z$ & $T$ & $U$ & $V$ & $E$ & $E$ \\ 
$E$ & $O$ & $X$ & $Y$ & $Z$ & $T$ & $U$ & $V$ & $E$%
\end{tabular}%
\end{equation*}

\textbf{2)} \textbf{Partially ordered case. }There is only two types of
partially ordered Wajsberg algebra with $8$ elements, up to an isomorphism
of ordered sets. Indeed, $\pi _{8}=2$. Let \newline
$\left( W_{1}=\{0,a,b,e\},\circ ,\overline{},e\right) $ and $\left(
W_{2}=\{0,1\},\cdot ,^{\prime },1\right) $ be two finite totally ordered
Wajsberg algebras. Using relation $\left( 3.1\right) $, on $W_{1}$ we have
that $\overline{b}=a$ and $\overline{a}=b$. We consider $W_{1}\times
W_{2}=\{\left( 0,0\right) ,\left( 0,1\right) ,\left( a,0\right) ,\left(
a,1\right) ,\left( b,0\right) ,\left( b,1\right) ,\left( e,0\right) ,\left(
e,1\right) \}=$\newline
$=\{O,X,Y,Z,T,U,V,E\}$. On $W_{1}\times W_{2}$ we obtain a Wajsberg algebra
structure by defining the multiplication as in relation $\left( 3.2\right) $%
, namely $\mathcal{W}_{11}^{8}=\left( W_{1}\times W_{2},\nabla
_{11}^{8}\right) $. The multiplication $\nabla _{11}^{8}$ is given in the
following table:%
\begin{equation}
\begin{tabular}{l|llllllll}
$\nabla _{11}^{8}$ & $O$ & $X$ & $Y$ & $Z$ & $T$ & $U$ & $V$ & $E$ \\ \hline
$O$ & $E$ & $E$ & $E$ & $E$ & $E$ & $E$ & $E$ & $E$ \\ 
$X$ & $V$ & $E$ & $V$ & $E$ & $V$ & $E$ & $V$ & $E$ \\ 
$Y$ & $U$ & $U$ & $E$ & $E$ & $E$ & $E$ & $E$ & $E$ \\ 
$Z$ & $T$ & $U$ & $V$ & $E$ & $V$ & $E$ & $V$ & $E$ \\ 
$T$ & $Z$ & $Z$ & $U$ & $U$ & $E$ & $E$ & $E$ & $E$ \\ 
$U$ & $Y$ & $Z$ & $T$ & $U$ & $T$ & $E$ & $V$ & $E$ \\ 
$V$ & $X$ & $X$ & $Z$ & $Z$ & $U$ & $U$ & $E$ & $E$ \\ 
$E$ & $O$ & $X$ & $Y$ & $Z$ & $T$ & $U$ & $V$ & $E$%
\end{tabular}%
.  \tag{4.3.}
\end{equation}%
(see [FV; 19], Example 4.13, relation $\left( 4.8\right) $)

In $\mathcal{W}_{11}^{8}$ we have that $O\leq X\leq Z\leq U\leq E$, $O\leq
Y\leq T\leq V\leq E$, $O\leq Y\leq Z\leq U\leq E,O\leq Y\leq T\leq U\leq E$
and the other elements can't be compared in this algebra. We denote this
order relation with $\leq _{11}^{8}$. All proper ideals $P_{1}=\{O,Y,T,V\}$, 
$P_{2}=\{O,X\}$ are prime ideals.\newline
We also have:\newline
-$\mathcal{W}_{11}^{8}/P_{1}=\{\overleftrightarrow{O},\overleftrightarrow{E}%
\}$ and $\mathcal{W}_{11}^{8}/P_{2}=\{\overline{O},\overline{Y},\overline{U},%
\overline{E}\}$, where $\overleftrightarrow{O}=\{O,Y,T,V\},%
\overleftrightarrow{E}=\{X,Z,U,E\},\overline{O}=\{O,X\},\overline{Y}%
=\{Y,Z\}, $\newline
$\overline{U}=\{U,T\},\overline{E}=\{V,E\}$;\newline
-$\mathcal{W}_{11}^{8}\simeq \mathcal{W}_{11}^{8}/P_{1}\times \mathcal{W}%
_{11}^{8}/P_{2}$, as in Remark 3.7, iv).

Now, we consider $W_{2}\times W_{1}=\{\left( 0,0\right) ,\left( 0,a\right)
,\left( 0,b\right) ,\left( 0,e\right) ,\left( 1,0\right) ,\left( 1,a\right)
,\left( 1,b\right) ,\left( 1,e\right) \}=$\newline
$=\{O,X,Y,Z,T,U,V,E\}$. On $W_{2}\times W_{1}$ we obtain a Wajsberg algebra
structure by defining the multiplication as in relation $\left( 3.2\right) $%
, namely $\mathcal{W}_{12}^{8}=\left( W_{2}\times W_{1},\nabla
_{12}^{8}\right) $. The multiplication $\nabla _{12}^{8}$ is given in
relation $\left( 4.8\right) $ from [FV;19], Example 4.13. We have that $%
O\leq X\leq Y\leq Z\leq E$, $O\leq X\leq Y\leq V\leq E$,\newline
$O\leq X\leq U\leq V\leq E$, $O\leq T\leq U\leq V\leq E$, and the other
elements can't be compared in this algebra. We denote this order relation
with $\leq _{12}^{8}.$ These two structures, $\mathcal{W}_{11}^{8}$ and $%
\mathcal{W}_{12}^{8}$, are isomorphic as Wajsberg algebras. The morphism is 
\begin{eqnarray*}
f_{12}^{8} &:&\left( W_{1}\times W_{2},\nabla _{11}^{8}\right) \rightarrow
\left( W_{2}\times W_{1},\nabla _{12}^{8}\right) , \\
f_{12}^{8}\left( X\right) &=&T,f_{12}^{8}\left( Y\right) =X,f_{12}^{8}\left(
Z\right) =U,f_{12}^{8}\left( T\right) =Y, \\
f_{12}^{8}\left( U\right) &=&V,f_{12}^{8}\left( V\right) =Z,f_{12}^{8}\left(
O\right) =O,f_{12}^{8}\left( E\right) =E.
\end{eqnarray*}%
All proper ideals $P_{1}=\{O,X,Y,Z\}$, $P_{2}=\{O,T\}$ are prime ideals.%
\newline
We also have:\newline
-$\mathcal{W}_{12}^{8}/P_{1}=\{\overleftrightarrow{O},\overleftrightarrow{E}%
\}$ and $\mathcal{W}_{12}^{8}/P_{2}=\{\overline{O},\overline{Y},\overline{U},%
\overline{E}\}$, where $\overleftrightarrow{O}=\{O,X,Y,Z\},%
\overleftrightarrow{E}=\{T,U,V,E\},\overline{O}=\{O,T\},\overline{U}%
=\{X,U\}, $\newline
$\overline{Y}=\{V,Y\},\overline{E}=\{Z,E\}$;\newline
-$\mathcal{W}_{12}^{8}\simeq \mathcal{W}_{12}^{8}/P_{1}\times \mathcal{W}%
_{12}^{8}/P_{2}$, as in Remark 3.7, iv).

If we take the map 
\begin{eqnarray*}
f_{13}^{8} &:&\left( W_{1}\times W_{2},\nabla _{11}^{8}\right) \rightarrow
\left( W_{2}\times W_{1},\nabla _{13}^{8}\right) , \\
f_{13}^{8}\left( X\right) &=&Y,f_{13}^{8}\left( Y\right) =U,f_{13}^{8}\left(
Z\right) =X,f_{13}^{8}\left( T\right) =V, \\
f_{13}^{8}\left( U\right) &=&Z,f_{13}^{8}\left( V\right) =T,f_{13}^{8}\left(
O\right) =O,f_{13}^{8}\left( E\right) =E,
\end{eqnarray*}%
we obtain the Wajsberg algebra $\mathcal{W}_{13}^{8}=\left( W_{2}\times
W_{1},\nabla _{13}^{8}\right) $. These two structures, $\mathcal{W}_{11}^{8}$
and $\mathcal{W}_{13}^{8}$, are isomorphic as ordered sets. In this algebra,
we have $O\leq Y\leq X\leq Z\leq E,O\leq U\leq V\leq T\leq E,$

$O\leq U\leq X\leq Z\leq E,O\leq U\leq V\leq Z\leq E$, and the other
elements can't be compared in this algebra. We denote this order relation
with $\leq _{13}^{8}$. The multiplication are given in the following table:%
\begin{equation*}
\begin{tabular}{l|llllllll}
$\nabla _{13}^{8}$ & $O$ & $X$ & $Y$ & $Z$ & $T$ & $U$ & $V$ & $E$ \\ \hline
$O$ & $E$ & $E$ & $E$ & $E$ & $E$ & $E$ & $E$ & $E$ \\ 
$X$ & $V$ & $E$ & $Z$ & $E$ & $T$ & $T$ & $T$ & $E$ \\ 
$Y$ & $T$ & $E$ & $E$ & $E$ & $T$ & $T$ & $T$ & $E$ \\ 
$Z$ & $U$ & $Z$ & $X$ & $E$ & $T$ & $V$ & $V$ & $E$ \\ 
$T$ & $Y$ & $X$ & $Y$ & $Z$ & $E$ & $X$ & $Z$ & $E$ \\ 
$U$ & $Z$ & $E$ & $Z$ & $E$ & $E$ & $E$ & $E$ & $E$ \\ 
$V$ & $X$ & $Z$ & $X$ & $E$ & $E$ & $Z$ & $E$ & $E$ \\ 
$E$ & $O$ & $X$ & $Y$ & $Z$ & $T$ & $U$ & $V$ & $E$%
\end{tabular}%
\end{equation*}

\begin{eqnarray*}
f_{13}^{8} &:&\left( W_{1}\times W_{2},\nabla _{11}^{8}\right) \rightarrow
\left( W_{2}\times W_{1},\nabla _{13}^{8}\right) , \\
f_{13}^{8}\left( X\right) &=&Y,f_{13}^{8}\left( Y\right) =U,f_{13}^{8}\left(
Z\right) =X,f_{13}^{8}\left( T\right) =V, \\
f_{13}^{8}\left( U\right) &=&Z,f_{13}^{8}\left( V\right) =T,f_{13}^{8}\left(
O\right) =O,f_{13}^{8}\left( E\right) =E.
\end{eqnarray*}

All proper ideals $P_{1}=\{O,U,V,T\}$, $P_{2}=\{O,Y\}$ are prime ideals.%
\newline
We also have:\newline
-$\mathcal{W}_{13}^{8}/P_{1}=\{\overleftrightarrow{O},\overleftrightarrow{E}%
\}$ and $\mathcal{W}_{13}^{8}/P_{2}=\{\overline{O},\overline{U},\overline{Z},%
\overline{E}\}$, where $\overleftrightarrow{O}=\{O,U,T,V\},%
\overleftrightarrow{E}=\{X,Z,Y,E\},\overline{O}=\{O,Y\},\overline{Y}%
=\{U,X\}, $\newline
$\overline{Z}=\{Z,V\},\overline{E}=\{V,E\}$;\newline
-$\mathcal{W}_{13}^{8}\simeq \mathcal{W}_{13}^{8}/P_{1}\times \mathcal{W}%
_{13}^{8}/P_{2}$, as in Remark 3.7, iv).

If we take $W_{2}\times W_{2}\times W_{2}=\{\left( 0,0,0\right) ,\left(
0,0,1\right) ,\left( 0,1,0\right) ,\left( 0,1,1\right) ,\left( 1,0,0\right)
, $\newline
$\left( 1,0,1\right) ,\left( 1,1,0\right) ,\left( 1,1,1\right)
\}=\{O,X,Y,Z,T,U,V,E\}$, on $W_{2}\times W_{2}\times W_{2}$ we obtain a
Wajsberg algebra structure by defining the multiplication as in relation $%
\left( 3.2\right) $, namely $\mathcal{W}_{21}^{8}=\left( W_{2}\times
W_{2}\times W_{2},\nabla _{21}^{8}\right) $. The multiplication $\nabla
_{21}^{8}$ is given in the following table:%
\begin{equation}
\begin{tabular}{l|llllllll}
$\nabla _{21}^{8}$ & $O$ & $X$ & $Y$ & $Z$ & $T$ & $U$ & $V$ & $E$ \\ \hline
$O$ & $E$ & $E$ & $E$ & $E$ & $E$ & $E$ & $E$ & $E$ \\ 
$X$ & $V$ & $E$ & $V$ & $E$ & $V$ & $E$ & $V$ & $E$ \\ 
$Y$ & $U$ & $U$ & $E$ & $E$ & $U$ & $U$ & $E$ & $E$ \\ 
$Z$ & $T$ & $U$ & $V$ & $E$ & $T$ & $U$ & $V$ & $E$ \\ 
$T$ & $Z$ & $Z$ & $Z$ & $Z$ & $E$ & $E$ & $E$ & $E$ \\ 
$U$ & $Y$ & $Z$ & $Y$ & $Z$ & $V$ & $E$ & $V$ & $E$ \\ 
$V$ & $X$ & $X$ & $Z$ & $Z$ & $U$ & $U$ & $E$ & $E$ \\ 
$E$ & $O$ & $X$ & $Y$ & $Z$ & $T$ & $U$ & $V$ & $E$%
\end{tabular}
\tag{4.4.}
\end{equation}%
(see [FV; 19], Example 4.13, relation $\left( 4.9\right) $).

We have that $X\leq Z,X\leq U,Y\leq Z,T\leq V,Y\leq V,T\leq U$ and the other
elements can't be compared in this algebra. We denote this order relation
with $\leq _{21}^{8}$. These two structures, $\mathcal{W}_{11}^{8}$ and $%
\mathcal{W}_{21}^{8}$, are not isomorphic as ordered sets, neither as
Wajsberg algebras. Here, we must remark that $W_{2}\times W_{2}\times
W_{2}\simeq \mathcal{W}_{21}^{8}$ and $W_{2}\times (W_{2}\times W_{2})$ $%
\simeq \mathcal{W}_{11}^{8}$. All proper ideals $P_{1}=\{O,X\}$, $%
P_{2}=\{O,Y\}$, $P_{3}=\{O,T\}$, $P_{4}=\{O,X,Y,Z\}$, $P_{5}=\{O,X,T,U\}$, $%
P_{6}=\{O,Y,T,V\}$ are also prime ideals.

We have that:\newline
- $\mathcal{W}_{21}^{8}/P_{1}=\{\overline{O},\overline{Y},\overline{T},%
\overline{E}\}$, where $\overline{O}=\{O,X\}$, $\overline{Y}=\{Y,Z\}$, $%
\overline{T}=\{T,U\},\overline{E}=\{V,E\}$;\newline
-$\mathcal{W}_{21}^{8}/P_{2}=\{\widehat{O},\widehat{X},\widehat{V},\widehat{E%
}\}$, where $\widehat{O}=\{O,Y\}$, $\widehat{X}=\{X,Z\}$, $\widehat{V}%
=\{V,U\}$, $\widehat{E}=\{T,E\}$.\newline
-$\mathcal{W}_{21}^{8}/P_{3}=\{\widetilde{O},\widetilde{X},\widetilde{Y},%
\widetilde{E}\}$, where $\widetilde{O}=\{O,T\}$, $\widetilde{X}=\{X,U\}$, $%
\widetilde{Y}=\{Y,V\}$, $\widetilde{E}=\{Z,E\}$;\newline
- $\mathcal{W}_{21}^{8}/P_{4}=\{\overleftrightarrow{O},\overleftrightarrow{E}%
\}$, where $\overleftrightarrow{O}=\{O,X,Y,Z\}$, $\overleftrightarrow{E}%
=\{U,V,T,E\}$; -$\mathcal{W}_{21}^{8}/P_{5}=\{\overleftarrow{O},%
\overleftarrow{E}\}$, where $\overleftarrow{O}=\{O,X,T,U\}$, $\overleftarrow{%
E}=\{Y,Z,V,E\}$;\newline
-$\mathcal{W}_{21}^{8}/P_{6}=\{\overrightarrow{O},\overrightarrow{E}\}$,
where $\overrightarrow{O}=\{O,Y,T,V\}$, $\overrightarrow{E}=\{X,Z,U,E\}$.

Therefore 
\begin{equation*}
\mathcal{W}_{11}^{8}\simeq \mathcal{W}_{21}^{8}/P_{1}\times \mathcal{W}%
_{21}^{8}/P_{6}\simeq \mathcal{W}_{21}^{8}/P_{2}\times \mathcal{W}%
_{21}^{8}/P_{5}\simeq \mathcal{W}_{21}^{8}/P_{3}\times \mathcal{W}%
_{21}^{8}/P_{4}
\end{equation*}
and%
\begin{equation*}
\mathcal{W}_{21}^{8}\simeq \mathcal{W}_{21}^{8}/P_{4}\times \mathcal{W}%
_{21}^{8}/P_{5}\times \mathcal{W}_{21}^{8}/P_{6},
\end{equation*}
as in Remark 3.7, iv).

If we consider the map

\begin{eqnarray*}
f_{22}^{8} &:&\left( W_{2}\times W_{2}\times W_{2},\nabla _{21}^{8}\right)
\rightarrow \left( W_{2}\times W_{2}\times W_{2},\nabla _{22}^{8}\right) , \\
f_{22}^{8}\left( X\right) &=&U,f_{22}^{8}\left( Y\right) =T,f_{22}^{8}\left(
Z\right) =V,f_{22}^{8}\left( T\right) =X, \\
f_{22}^{8}\left( U\right) &=&Z,f_{22}^{8}\left( V\right) =Y,f_{22}^{8}\left(
O\right) =O,f_{22}^{8}\left( E\right) =E,
\end{eqnarray*}%
we obtain the Wajsberg algebra $\mathcal{W}_{22}^{8}=\left( W_{2}\times
W_{2}\times W_{2},\nabla _{22}^{8}\right) $ given in the below table:

\begin{equation*}
\begin{tabular}{l|llllllll}
$\nabla _{22}^{8}$ & $O$ & $X$ & $Y$ & $Z$ & $T$ & $U$ & $V$ & $E$ \\ \hline
$O$ & $E$ & $E$ & $E$ & $E$ & $E$ & $E$ & $E$ & $E$ \\ 
$X$ & $V$ & $E$ & $E$ & $E$ & $V$ & $V$ & $V$ & $E$ \\ 
$Y$ & $U$ & $Z$ & $E$ & $Z$ & $V$ & $U$ & $V$ & $E$ \\ 
$Z$ & $T$ & $Y$ & $Y$ & $E$ & $T$ & $V$ & $V$ & $E$ \\ 
$T$ & $Z$ & $Z$ & $E$ & $Z$ & $E$ & $Z$ & $E$ & $E$ \\ 
$U$ & $Y$ & $Y$ & $Y$ & $E$ & $Y$ & $E$ & $E$ & $E$ \\ 
$V$ & $X$ & $X$ & $Y$ & $Z$ & $Y$ & $Z$ & $E$ & $E$ \\ 
$E$ & $O$ & $X$ & $Y$ & $Z$ & $T$ & $U$ & $V$ & $E$%
\end{tabular}%
\end{equation*}

These two structures, $\mathcal{W}_{21}^{8}$ and $\mathcal{W}_{22}^{8}$, are
isomorphic as Wajsberg algebras. We have that $U\leq V,U\leq Z,T\leq V,X\leq
Y,T\leq Y,X\leq Z$ and the other elements can't be compared in this algebra.
We denote this order relation with $\leq _{22}^{8}$.

All proper ideals $P_{1}=\{O,U\}$, $P_{2}=\{O,T\}$, $P_{3}=\{O,X\}$, $%
P_{4}=\{O,U,T,V\}$, $P_{5}=\{O,U,X,Z\}$, $P_{6}=\{O,T,X,Y\}$ are also prime
ideals.

We obtain that:\newline
- $\mathcal{W}_{22}^{8}/P_{1}=\{\overline{O},\overline{T},\overline{X},%
\overline{E}\}$, where $\overline{O}=\{O,U\}$, $\overline{T}=\{T,V\}$, $%
\overline{X}=\{X,Z\},\overline{E}=\{Y,E\}$;\newline
-$\mathcal{W}_{22}^{8}/P_{2}=\{\widehat{O},\widehat{U},\widehat{Z},\widehat{E%
}\}$, where $\widehat{O}=\{O,T\}$, $\widehat{U}=\{U,V\}$, $\widehat{Z}%
=\{Z,Y\}$, $\widehat{E}=\{X,E\}$.\newline
-$\mathcal{W}_{22}^{8}/P_{3}=\{\widetilde{O},\widetilde{Z},\widetilde{Y},%
\widetilde{E}\}$, where $\widetilde{O}=\{O,X\}$, $\widetilde{Z}=\{Z,U\}$, $%
\widetilde{Y}=\{Y,T\}$, $\widetilde{E}=\{V,E\}$;\newline
- $\mathcal{W}_{22}^{8}/P_{4}=\{\overleftrightarrow{O},\overleftrightarrow{E}%
\}$, where $\overleftrightarrow{O}=\{O,U,T,V\}$, $\overleftrightarrow{E}%
=\{X,Y,Z,E\}$; -$\mathcal{W}_{22}^{8}/P_{5}=\{\overleftarrow{O},%
\overleftarrow{E}\}$, where $\overleftarrow{O}=\{O,U,X,Z\}$, $\overleftarrow{%
E}=\{Y,V,T,E\}$;\newline
-$\mathcal{W}_{22}^{8}/P_{6}=\{\overrightarrow{O},\overrightarrow{E}\}$,
where $\overrightarrow{O}=\{O,T,X,Y\}$, $\overrightarrow{E}=\{V,Z,U,E\}$.

Therefore 
\begin{equation*}
\mathcal{W}_{11}^{8}\simeq \mathcal{W}_{22}^{8}/P_{1}\times \mathcal{W}%
_{22}^{8}/P_{6}\simeq \mathcal{W}_{22}^{8}/P_{2}\times \mathcal{W}%
_{22}^{8}/P_{5}\simeq \mathcal{W}_{22}^{8}/P_{3}\times \mathcal{W}%
_{22}^{8}/P_{4}
\end{equation*}%
and%
\begin{equation*}
\mathcal{W}_{21}^{8}\simeq \mathcal{W}_{22}^{8}/P_{4}\times \mathcal{W}%
_{22}^{8}/P_{5}\times \mathcal{W}_{22}^{8}/P_{6},
\end{equation*}%
as in Remark 3.7, iv).

If we take the map

\begin{eqnarray*}
f_{23}^{8} &:&\left( W_{2}\times W_{2}\times W_{2},\nabla _{21}^{8}\right)
\rightarrow \left( W_{2}\times W_{2}\times W_{2},\nabla _{23}^{8}\right) , \\
f_{23}^{8}\left( X\right) &=&Z,f_{23}^{8}\left( Y\right) =X,f_{23}^{8}\left(
Z\right) =V,f_{23}^{8}\left( T\right) =U, \\
f_{23}^{8}\left( U\right) &=&T,f_{23}^{8}\left( V\right) =Y,f_{23}^{8}\left(
O\right) =O,f_{23}^{8}\left( E\right) =E,
\end{eqnarray*}%
we obtain the Wajsberg algebra $\mathcal{W}_{23}^{8}=\left( W_{2}\times
W_{2}\times W_{2},\nabla _{23}^{8}\right) $ given in the below table:%
\begin{equation*}
\begin{tabular}{l|llllllll}
$\nabla _{23}^{8}$ & $O$ & $X$ & $Y$ & $Z$ & $T$ & $U$ & $V$ & $E$ \\ \hline
$O$ & $E$ & $E$ & $E$ & $E$ & $E$ & $E$ & $E$ & $E$ \\ 
$X$ & $T$ & $E$ & $E$ & $T$ & $T$ & $T$ & $E$ & $E$ \\ 
$Y$ & $Z$ & $V$ & $E$ & $Z$ & $T$ & $T$ & $V$ & $E$ \\ 
$Z$ & $Y$ & $Y$ & $Y$ & $E$ & $E$ & $Y$ & $E$ & $E$ \\ 
$T$ & $X$ & $X$ & $Y$ & $V$ & $E$ & $Y$ & $V$ & $E$ \\ 
$U$ & $V$ & $V$ & $E$ & $V$ & $E$ & $E$ & $V$ & $E$ \\ 
$V$ & $U$ & $Y$ & $Y$ & $T$ & $T$ & $U$ & $E$ & $E$ \\ 
$E$ & $O$ & $X$ & $Y$ & $Z$ & $T$ & $U$ & $V$ & $E$%
\end{tabular}%
.
\end{equation*}%
These two structures, $\mathcal{W}_{21}^{8}$ and $\mathcal{W}_{23}^{8}$, are
isomorphic as ordered sets.We have that $Z\leq V,Z\leq T,X\leq V,U\leq
Y,X\leq Y,U\leq T$ and the other elements can't be compared in this algebra.
We denote this order relation with $\leq _{23}^{8}$.

All proper ideals $P_{1}=\{O,Z\}$, $P_{2}=\{O,X\}$, $P_{3}=\{O,U\}$, $%
P_{4}=\{O,Z,X,V\}$, $P_{5}=\{O,Z,U,T\}$, $P_{6}=\{O,X,U,Y\}$ are also prime
ideals.

We obtain that:\newline
- $\mathcal{W}_{23}^{8}/P_{1}=\{\overline{O},\overline{X},\overline{T},%
\overline{E}\}$, where $\overline{O}=\{O,Z\}$, $\overline{X}=\{X,V\}$, $%
\overline{T}=\{T,U\},\overline{E}=\{Y,E\}$;\newline
-$\mathcal{W}_{23}^{8}/P_{2}=\{\widehat{O},\widehat{V},\widehat{Y},\widehat{E%
}\}$, where $\widehat{O}=\{O,X\}$, $\widehat{V}=\{V,Z\}$, $\widehat{Y}%
=\{Y,T\}$, $\widehat{E}=\{U,E\}$.\newline
-$\mathcal{W}_{23}^{8}/P_{3}=\{\widetilde{O},\widetilde{Z},\widetilde{Y},%
\widetilde{E}\}$, where $\widetilde{O}=\{O,U\}$, $\widetilde{}=\{Z,T\}$, $%
\widetilde{Y}=\{Y,X\}$, $\widetilde{E}=\{V,E\}$;\newline
- $\mathcal{W}_{23}^{8}/P_{4}=\{\overleftrightarrow{O},\overleftrightarrow{E}%
\}$, where $\overleftrightarrow{O}=\{O,Z,X,V\}$, $\overleftrightarrow{E}%
=\{U,Y,T,E\}$; -$\mathcal{W}_{23}^{8}/P_{5}=\{\overleftarrow{O},%
\overleftarrow{E}\}$, where $\overleftarrow{O}=\{O,Z,U,T\}$, $\overleftarrow{%
E}=\{Y,X,V,E\}$;\newline
-$\mathcal{W}_{23}^{8}/P_{6}=\{\overrightarrow{O},\overrightarrow{E}\}$,
where $\overrightarrow{O}=\{O,X,U,Y\}$, $\overrightarrow{E}=\{V,Z,T,E\}$.

Therefore 
\begin{equation*}
\mathcal{W}_{11}^{8}\simeq \mathcal{W}_{23}^{8}/P_{1}\times \mathcal{W}%
_{23}^{8}/P_{6}\simeq \mathcal{W}_{23}^{8}/P_{2}\times \mathcal{W}%
_{23}^{8}/P_{5}\simeq \mathcal{W}_{23}^{8}/P_{3}\times \mathcal{W}%
_{23}^{8}/P_{4}
\end{equation*}%
and%
\begin{equation*}
\mathcal{W}_{21}^{8}\simeq \mathcal{W}_{23}^{8}/P_{4}\times \mathcal{W}%
_{23}^{8}/P_{5}\times \mathcal{W}_{23}^{8}/P_{6},
\end{equation*}%
as in Remark 3.7, iv).

From the above, we have that there are only three types of nonisomorphic
Wajsberg algebras of order $8$. Using connections between MV-algebras and
Wajsberg algebras, we obtain that there are two nonisomorphic MV-algebras of
order $8$. The same number was found in [BV; 10], Table 8, by using other
method. Now, we count the number of all Wajsberg algebras of order $8$. For
this purpose, we must find all isomofirsms of ordered sets $f:$ $\mathcal{W}%
_{11}^{8}\rightarrow \mathcal{W}_{11}^{8}$, such that $f\left( O\right)
=O,f\left( E\right) =E$ and all isomofirsms of ordered sets $f:$ $\mathcal{W}%
_{21}^{8}\rightarrow \mathcal{W}_{21}^{8}$, such that $f\left( O\right)
=O,f\left( E\right) =E$. Therefore, there are $2\times 6!=1440$ isomorphisms
and, in turn, $1440$ partially ordered Wajsberg algebras. In total, there
are $1441$ Wajsberg algebras of order $8$.\bigskip 

\textbf{4.4.} \textbf{Wajsberg algebras of order }$9\bigskip $

1) \textbf{Totally ordered case. }Let $W=\{O\leq X\leq Y\leq Z\leq T\leq
U\leq S\leq V\leq E\}$ be a totally ordered set. On $W$ we define a
multiplication as in relation $\left( 3\mathbf{.}1\right) $. We have $%
\overline{X}=V$, $\overline{Y}=S$, $\overline{Z}=U,\overline{T}=T$.
Therefore the algebra $W$ has the following multiplication table:%
\begin{equation*}
\begin{tabular}{l|lllllllll}
$\nabla _{0}^{9}$ & $O$ & $X$ & $Y$ & $Z$ & $T$ & $U$ & $S$ & $V$ & $E$ \\ 
\hline
$O$ & $E$ & $E$ & $E$ & $E$ & $E$ & $E$ & $E$ & $E$ & $E$ \\ 
$X$ & $V$ & $E$ & $E$ & $E$ & $E$ & $E$ & $E$ & $E$ & $E$ \\ 
$Y$ & $S$ & $V$ & $E$ & $E$ & $E$ & $E$ & $E$ & $E$ & $E$ \\ 
$Z$ & $U$ & $S$ & $V$ & $E$ & $E$ & $E$ & $E$ & $E$ & $E$ \\ 
$T$ & $T$ & $U$ & $S$ & $V$ & $E$ & $E$ & $V$ & $E$ & $E$ \\ 
$U$ & $Z$ & $T$ & $U$ & $S$ & $V$ & $E$ & $S$ & $V$ & $E$ \\ 
$S$ & $Y$ & $Z$ & $T$ & $U$ & $S$ & $V$ & $E$ & $E$ & $E$ \\ 
$V$ & $X$ & $Y$ & $Z$ & $T$ & $U$ & $S$ & $V$ & $E$ & $E$ \\ 
$E$ & $O$ & $X$ & $Y$ & $Z$ & $T$ & $U$ & $S$ & $V$ & $E$%
\end{tabular}%
\end{equation*}

\textbf{2)} \textbf{Partially ordered case. }There is only one type of
partially ordered Wajsberg algebra with $9$ elements, up to an isomorphism
of ordered sets. Indeed, $\pi _{9}=1\left( W_{1}=\{0,a,e\},\circ ,\overline{}%
,e\right) $ and $\left( W_{2}=\{0,b,1\},\cdot ,^{\prime },1\right) .$We have 
$\overline{a}=a$ and $b^{\prime }=b.$ We have \newline
$W_{1}\times W_{2}=\{\left( 0,0\right) ,\left( 0,b\right) ,\left( 0,1\right)
,\left( a,0\right) ,\left( a,b\right) ,\left( a,1\right) ,\left( e,0\right)
,\left( e,b\right) ,\left( e,1\right) \}=$\newline
$=\{O,X,Y,Z,T,U,S,V,E\}.$

\begin{equation*}
\begin{tabular}{l|lllllllll}
$\nabla _{11}^{9}$ & $O$ & $X$ & $Y$ & $Z$ & $T$ & $U$ & $S$ & $V$ & $E$ \\ 
\hline
$O$ & $E$ & $E$ & $E$ & $E$ & $E$ & $E$ & $E$ & $E$ & $E$ \\ 
$X$ & $V$ & $E$ & $E$ & $V$ & $E$ & $E$ & $V$ & $E$ & $E$ \\ 
$Y$ & $S$ & $V$ & $E$ & $S$ & $V$ & $E$ & $S$ & $V$ & $E$ \\ 
$Z$ & $U$ & $U$ & $U$ & $E$ & $E$ & $E$ & $E$ & $E$ & $E$ \\ 
$T$ & $T$ & $U$ & $U$ & $V$ & $E$ & $E$ & $V$ & $E$ & $E$ \\ 
$U$ & $Z$ & $T$ & $U$ & $S$ & $V$ & $E$ & $S$ & $V$ & $E$ \\ 
$S$ & $Y$ & $X$ & $Y$ & $U$ & $U$ & $U$ & $E$ & $E$ & $E$ \\ 
$V$ & $X$ & $Y$ & $Y$ & $T$ & $U$ & $U$ & $V$ & $E$ & $E$ \\ 
$E$ & $O$ & $X$ & $Y$ & $Z$ & $T$ & $U$ & $S$ & $V$ & $E$%
\end{tabular}%
\end{equation*}

We obtain a Wajsberg algebra structure, namely $\mathcal{W}_{11}^{9}=\left(
W_{1}\times W_{2},\nabla _{11}^{9}\right) $. We remark that $O\leq X\leq
Y\leq U,O\leq X\leq T\leq U,O\leq X\leq T\leq V,$\newline
$O\leq Z\leq T\leq U,O\leq Z\leq T\leq V,O\leq Z\leq S\leq V$ and the other
elements can't be compared. We denote this order relation with $\leq
_{11}^{9}$. All proper ideals $P_{1}=\{O,X,Y\}$, $P_{2}=\{O,Z,S\}$, are also
prime ideals.

We also have:\newline
-$\mathcal{W}_{11}^{9}/P_{1}=\{\overleftrightarrow{O},\overleftrightarrow{U},%
\overleftrightarrow{E}\}$ and $\mathcal{W}_{11}^{9}/P_{2}=\{\overline{O},%
\overline{T},\overline{E}\}$, where $\overleftrightarrow{O}=\{O,X,Y\},%
\overleftrightarrow{U}=\{U,T,Z\},\overleftrightarrow{E}=\{S,V,E\},$\newline
$\overline{O}=\{O,Z,S\},\overline{T}=\{T,V,X\},\overline{E}=\{Y,U,E\}$;%
\newline
-$\mathcal{W}_{11}^{9}\simeq \mathcal{W}_{11}^{9}/P_{1}\times \mathcal{W}%
_{11}^{9}/P_{2}$, as in Remark 3.7, iv).

From the above, we have that there are only two types of nonisomorphic
Wajsberg algebras of order $9$. Using connections between MV-algebras and
Wajsberg algebras, we obtain that there are two nonisomorphic MV-algebras of
order $9$. The same number was found in [BV; 10], Table 8, by using other
method. Now, we count the number of Wajsberg algebras of order $9$. For this
purpose, we must find all isomofirsms of ordered sets $f:$ $\mathcal{W}%
_{11}^{9}\rightarrow \mathcal{W}_{11}^{9}$, such that $f\left( O\right)
=O,f\left( E\right) =E$, such that $f\left( O\right) =O,f\left( E\right) =E$%
. Therefore, there are $7!=5040$ isomorphisms and, in turn, $5040$ partially
ordered Wajsberg algebras. In total, there are $5041$ Wajsberg algebras of
order $9$.\medskip 

\textbf{Remark 4.1. }For $n\in \{2,3,5,7\}$, the Wajsberg algebras of order $%
n$ are totally ordered, since $n$ is a prime number. The multiplication
tables are easy to obtain.\medskip

From the above, we obtain the following proposition.\medskip

\textbf{Proposition 4.2.} \textit{With the above notations we have:}

1) \textit{The number of nonisomorphic Wajsberg algebras types of order} $n$ 
\textit{is} $\pi _{n}.$

2) \textit{The total number of Wajsberg algebras of order} $n$ \textit{is} $%
\pi _{n}\left( n-2\right) !+1$\textit{, for} $n\geq 4$\textit{,} $n$ \textit{%
not a prime number.} $\Box $

\begin{equation*}
\end{equation*}

\textbf{Conclusions\bigskip }

Starting from some results obtained in [FV;19], in this paper we provided a
Representation Theorem for finite Wajsberg algebras. Using this theorem, we
describe all finite Wajsberg algebras of order $n\leq 9$.

Wajsberg algebras have many applications in various domains and their study
can provide us new and interesting properties of them.

\begin{equation*}
\end{equation*}

\textbf{References}%
\begin{equation*}
\end{equation*}

[BV; 10] Belohlavek, R., Vilem Vychodil, V., \textit{Residuated Lattices of
Size} $\leq 12$, Order, 27(2010), 147-161.

[COM; 00] Cignoli, R. L. O, D'Ottaviano, I. M. L., Mundici, D., \textit{%
Algebraic foundations of many-valued reasoning}, Trends in Logic, Studia
Logica Library, Dordrecht, Kluwer Academic Publishers, 7(2000).

[CT; 96] Cignoli, R., Torell, A., T., \textit{Boolean Products of
MV-Algebras: Hypernormal MV-Algebras}, J Math Anal Appl (199)(1996), 637-653.

[CHA; 58] Chang, C.C.,\textit{\ Algebraic analysis of many-valued logic},
Trans. Amer. Math. Soc. 88(1958), 467-490.

[CHA; 59] Chang, C.C.,\textit{\ A new proof of the completeness of the
Lukasiewicz}, Trans. Amer. Math. Soc. 93(1959), 74-90.

[Di; 38] Dilworth, R.P., \textit{Abstract residuation over lattices}, Bull.
Am. Math. Soc. 44(1938), 262--268.

[FV; 19] Flaut, C, Vasile R., \textit{Wajsberg algebras arising from binary
block codes}, https://arxiv.org/pdf/1904.07169.pdf, 2019

[FRT; 84] Font, J., M., Rodriguez, A., J., Torrens, A., \textit{Wajsberg
Algebras}, Stochastica, 8(1)(1984), 5-30.

[GA; 90] Gaitan, H., \textit{About quasivarieties of p-algebras and Wajsberg
algebras}, 1990, Retrospective Theses and Dissertations, 9440,
https://lib.dr.iastate.edu/rtd/9440

[HR; 99] H\"{o}hle, U., Rodabaugh, S., E., \textit{Mathematics of Fuzzy
Sets: Logic, Topology and Measure Theory}, Springer Science and Business
Me\qquad qdia, LLC, 1999.

[Pi; 07] Piciu, D., \textit{Algebras of Fuzzy Logic}, Editura Universitaria,
Craiova, 2007.

[Tu; 99] Turunen, E.,\textit{\ Mathematics Behind Fuzzy Logic}, \textit{%
Advances in Soft Computing}, Heidelberg Physics, Verlag, 1999.

[WD; 39] Ward, M., Dilworth, R.P., \textit{Residuated lattices}, Trans. Am.
Math. Soc. 45(1939), 335--354.%
\begin{equation*}
\end{equation*}

Cristina Flaut

{\small Faculty of Mathematics and Computer Science, Ovidius University,}

{\small Bd. Mamaia 124, 900527, CONSTANTA, ROMANIA}

{\small \ http://www.univ-ovidius.ro/math/}

{\small e-mail: cflaut@univ-ovidius.ro; cristina\_flaut@yahoo.com}%
\begin{equation*}
\end{equation*}

\v{S}\'{a}rka Ho\v{s}kov\'{a}-Mayerov\'{a}

{\small Department of Mathematics and Physics,}

{\small University of Defence, Brno, Czech Republic}

{\small e-mail: sarka.mayerova@unob.cz}%
\begin{equation*}
\end{equation*}

Arsham Borumand Saeid

{\small Department of Pure Mathematics,}

{\small Faculty of Mathematics and Computer,}

{\small Shahid Bahonar University of Kerman, Kerman, Iran }

{\small e-mail: arsham@uk.ac.ir}%
\begin{equation*}
\end{equation*}%
\qquad \qquad

Radu Vasile,

{\small PhD student at Doctoral School of Mathematics,}

{\small Ovidius University of Constan\c{t}a, Rom\^{a}nia}

{\small rvasile@gmail.com}

\end{document}